\documentclass{article}
\usepackage[all]{xy}
\usepackage{graphicx}
\usepackage[margin=2.95cm]{geometry}
\usepackage{amssymb}
\usepackage{amsmath}
\usepackage{comment}
\usepackage{latexsym} 
\newtheorem{thm}{Theorem}[section]

\newtheorem{prop}[thm]{Proposition}

\newtheorem{rmk}[thm]{Remark}

\newcommand{\pf}{\noindent{\bf Proof~~}\ }

\newcommand{\reals}{{\mathbb R}}

\newcommand{\xbar}{{\overline{X}}}
\newcommand{\ybar}{{\overline{Y}}}

\newcommand{\qed}{\begin{flushright} $\Box$\ \ \ \ \ \end{flushright}}

\newcommand{\from}{\leftarrow}

\newcommand{\lrel}{\mathbf{LREL}}

\newcommand{\clrel}{\mathbf{CLREL}}
\newcommand{\ilrel}{\mathbf{ILREL}}
\newcommand{\slrel}{\mathbf{SLREL}}

\newcommand{\arrows}{\,\lower1pt\hbox{$\longrightarrow$}\hskip-.24in\raise2pt
             \hbox{$\longrightarrow$}\,}

\newcommand{\zero}{\mathbf{0}}

\usepackage{url}

\title{Decomposition of (co)isotropic relations}
\author{Jonathan Lorand\thanks{Research partially supported by a scholarship of the Anna \& Hans K\"agi Foundation. Current address: Institute of Mathematics, University of Zurich, Zurich, Switzerland.}\\
Department of Mathematics\\
ETH Zurich\\
Zurich, Switzerland\\
\and 
Alan Weinstein\thanks{Research partially supported by UC Berkeley Committee on Research
\newline\mbox{~~~~}MSC2010 Subject Classification Number: 
 18B10 (Primary), 53D99 (Secondary)
\newline \mbox{~~~~}Keywords: linear relation, duality, symplectic  vector space,
coisotropic relation}\\
Department of Mathematics\\University of California\\
Berkeley, CA 94720-3840 USA}
\date{}
\begin{document}
\maketitle

\begin{abstract}
We identify thirteen isomorphism classes of indecomposable coisotropic relations between Poisson vector spaces and show that every coisotropic relation between finite-dimensional Poisson vector spaces
may be decomposed as a direct sum of multiples of these indecomposables.   We also find a list of thirteen invariants, each of which is the dimension of a space constructed from the relation, such that the 
13-vector of multiplicities and the 13-vector of invariants are related by an invertible matrix over $\mathbb Z$.  

It turns out to be simpler to do the analysis above for isotropic relations between presymplectic vector spaces.   The coisotropic/Poisson case then follows by a simple duality argument.
\end{abstract}

\begin{center}
\emph{Dedicated to the memory of Louis Boutet de Monvel}
\end{center}

\section{Introduction}

The geometric quantization program in symplectic geometry and its implementation in the form of Fourier integral operators in microlocal analysis suggest that the natural symplectic analogues, or classical limits, of the linear operators between Hilbert spaces in quantum mechanics are the canonical relations between symplectic manifolds.  The relations, which are lagrangian submanifolds of products of the form $(X,\omega_X)\times
(Y,-\omega_Y),$  include the graphs of symplectic diffeomorphisms (the classical canonical transformations, which correspond to unitary operators via quantization) and much more.  In Poisson geometry, the graphs of Poisson maps are coisotropic, which makes coisotropic relations the natural ``maps'' to study in the Poisson setting.  Finally, duality suggests that for presymplectic manifolds, carrying closed forms which may be degenerate, the most important relations are the isotropic ones.

In each of these situations, it is tempting to treat the relevant relations as morphisms in a category, but they generally behave badly under arbitrary compositions (like the unbounded operators in quantum mechanics).  For instance, the set theoretic composition of two smooth relations can fail to be smooth.  The case of linear relations is somewhat better, since the classes of relations of interest are closed under composition, but even here there is a serious problem: the composition can depend discontinuously on the relations being composed.  Some approaches to this problem were investigated in 
 \cite{li-we:selective}  and \cite{we:(co)isotropic}.

In particular, in \cite{we:(co)isotropic}, the second author studied categories obtained from various categories of linear relations by ``regularizing" compositions which are degenerate in certain senses.   Each such category is an extension of the original one by a submonoid of the pairs of nonnegative integers. 
In the present paper, we analyze the individual relations themselves, showing that each one can be decomposed as a direct sum of relations of very simple types.   

The categories in which we work here are closely related to those in 
\cite{we:(co)isotropic}, namely $\lrel$, $\ilrel$, $\clrel$ and $\slrel$, where the first category is that of all linear relations between (finite dimensional) vector spaces, and the initial letter of each of the remaining three identifies the category of isotropic, coisotropic, or symplectic (\text{i.e.} lagrangian) relations between symplectic vector spaces.   
In fact, as was already noted in \cite{lo-we:coisotropic}, it is natural to consider Poisson vector spaces as the target and source objects for linear coisotropic relations (such as the graphs of linear Poisson maps), while the targets and sources in the isotropic case can be arbitrary presymplectic (\text{i.e.} carrying a possibly degenerate bilinear skew-symmetric form) vector spaces.   These two categories (and their double categories of commutative squares) are equivalent to one another via duality, and it will turn out to be convenient to prove our main results in the isotropic case and then transfer them to the coisotropic situation by duality.

In each case, we will study morphisms up to the equivalences given by pairs of linear isomorphisms between sources and targets taking one relation to another.   (These equivalences are the invertible 2-morphisms for a double category structure on each of our categories.)

Each of our categories is rigid monoidal, with direct sum $X \oplus Y$ as monoidal product.  
For clarity, we will use the symbol $\times$ and the term ``product"  when dealing with the direct sum of the target and the source space of a relation; we reserve the symbol $\oplus$ and the term ``direct sum" for decompositions \emph{within} the respective target and source spaces of a relation or for decompositions of relations themselves and their ambient product spaces.   The monoidal dual $\xbar$ is  the same vector space with the Poisson or presymplectic structure multiplied by -1. The unit object  in each case is the vector space  $\zero$ consisting of the zero element of the ground field $\reals$ (which may be replaced without any other changes by any field  of characteristic $\neq 2$).   

Our classification of morphisms will consist of a description of thirteen indecomposable morphisms and a proof that each morphism is an essentially unique orthogonal direct sum of multiples of these.
We will present a list of thirteen invariants of a relation, each of which is the dimension of a space constructed from the relation itself and the source and target Poisson or presymplectic spaces, and will show that the 13-vector of multiplicities is related to the 13-vector of invariants by multiplication by an invertible matrix over $\mathbb Z$. It follows that two morphisms are isomorphic (equivalent) if, and only if, their corresponding multiplicities are equal. 

Towber \cite{to:linear} carried out the classification for the case of $\lrel$, which may be thought of as the spaces with zero Poisson or presymplectic structure.
(His result is mentioned in passing as a remark after his solution of the much more substantial classification problem for endomorphisms up to conjugacy.) There are five isomorphism classes of indecomposable relations, represented by:
\begin{itemize}
\item  The identity relation $\reals \from \reals$;
\item  The zero relations $\reals \from \zero$ and $\zero \from \reals$;
\item  The relations $\reals \from \zero$ and $\zero \from \reals$ given by
 $\reals \times \zero$ and $\zero \times \reals$ respectively.
\end{itemize}
Any linear relation $X \from Y$ is isomorphic to a unique (up to order) direct sum of copies of these; the isomorphism itself is not unique; even the decompositions of $X$ and $Y$ are nonunique if the decomposables occur with multiplicities.   A suitable 5-vector of invariants is given by the dimensions of the source, target, relation, null space, and null space of the transpose.

In the symplectic (lagrangian) case, there are only three isomorphism classes of indecomposables, represented by:
\begin{itemize}
\item  The identity relation $\reals^2 \from \reals^2$;
\item The relations $\reals^2 \from \zero$ and $\zero \from \reals^2$ given by  $(\reals,0)\times \zero$
and $\zero \times (\reals,0)$ respectively.
\end{itemize}
This result follows immediately from the well-known fact\footnote{See, for instance, \cite{be-tu:relazioni}, Proposizioni 4.4 \& 4.5.} that every canonical relation is the direct sum of (the graph of) a symplectomorphism and a product of lagrangian subspaces in the target and source.

For isotropic and coisotropic relations between {\em symplectic} vector spaces, there are six isomorphism classes.  In the isotropic case, we will show that, in addition to the three lagrangian relations above, we have three further indecomposable types:
\begin{itemize}
\item  The zero relations $\reals^2 \from \zero$ and $\zero \from \reals^2$;
\item  The relation $\reals^2 \from \reals^2$ given in coordinates $(q_1,p_1,q_2,p_2)$ by the equations
$q_1 = q_2 =0$ and $p_1=p_2$.
\end{itemize}  

For the coisotropic case, we have merely to take symplectic orthogonals.  The three types beyond the lagrangian case are now:

\begin{itemize}
\item  The relations $\reals^2 \from \zero$ and $\zero \from \reals^2$ given by
 $\reals^2 \times \zero$ and $\zero \times \reals^2$ respectively.
\item  The relation $\reals^2 \from \reals^2$ given in coordinates $(q_1,p_1,q_2,p_2)$ by the equation
$p_1=p_2$.
\end{itemize}  

In the general case of isotropic relations between presymplectic spaces, one finds, in addition to the six isotropic relations between symplectic spaces listed above, the following five indecomposable types between spaces with the zero presymplectic structure:
 \begin{itemize}
\item The identity relation $\mathbb{R} \from \mathbb{R}$, 
\item The relations $\mathbb{R} \from \zero$ and $\zero \from \mathbb{R}$ given by $\mathbb{R} \times \zero$ and $\zero \times \mathbb{R}$ respectively,
\item The zero relations $\mathbb{R} \from \zero$ and $\zero \from \mathbb{R}$;
\end{itemize}
and the following two indecomposable types between a symplectic space and a zero-presymplectic space:
\begin{itemize}
\item The relations $\mathbb{R}^2 \from \mathbb{R}$ and $\mathbb{R} \from \mathbb{R}^2$ given by the natural isomorphisms $(\mathbb{R}, 0) \from \mathbb{R}$ and $\mathbb{R} \from (\mathbb{R},0)$ respectively. 
\end{itemize}

Note that the indecomposables where both target and source space carry the zero form are simply the indecomposables in $\lrel$, since every linear relation is isotropic with respect to the zero form. 

\bigskip
For coisotropic relations between Poisson vector spaces, it suffices to take the annihilators (in the dual product space) of the thirteen indecomposable types of isotropic relations between presymplectic spaces. A complete list 
for the Poisson case is given in Theorem \ref{thm-indecomp}.

The proofs and results in this paper are very similar to those for 
pairs of (co)isotropic subspaces in \cite{lo-we:coisotropic}.    We wonder whether these two decomposition theories might be examples of a general theory of representations of quivers by relations, which would be closely related to work of Sergeichuk \cite{se:classification systems} on representations of partially directed graphs by mappings and bilinear forms.   Here, as in \cite{lo-we:coisotropic}, we give constructions and descriptions of the decompositions which use nothing more than linear algebra.


\section{ $\text{Isotropic} = \text{cartesian} \oplus \text{biinjective} $}

As was the case for the classifications of subspace pairs in \cite{lo-we:coisotropic}, it is simplest to start with the case of isotropic relations between presymplectic spaces and pass by duality to the coisotropic relations between Poisson spaces.   

If $A$ is any subspace of a  space $X$ with presymplectic form $\omega$, we will denote by $A^\perp$ its presymplectic \textbf{orthogona}l $\{x\in X \mid\omega (x,a) = 0, ~\forall a\in A\}$.  
Thus, $A$ is isotropic when $A \subseteq A^\perp$.
The orthogonal  $X^\perp$ of $X$  itself is called the \textbf{radical} of $X$ and will be denoted by $R_X$.

For any  linear relation $X \overset{f}{\longleftarrow} Y$, its \textbf{kernel} $0f \subseteq Y$ is the subspace
$\{y \in Y \mid (0,y) \in f\}$, and its \textbf{indeterminacy} $f0 \subseteq X$ is $\{x\in X\mid(x,0) \in f\}$, the kernel of the transposed relation.   They are isotropic when $f$ is isotropic.   $f$ is \textbf{biinjective} when its kernel and indeterminacy are both zero and \textbf{cartesian} when it is the cartesian product of its image $fY$ and domain $Xf$.  (In the Poisson case, the image and domain are coisotropic when $f$ is coisotropic.) 

The following proposition is a generalization of the decomposition result for lagrangian relations cited in the introduction.

\begin{prop}
\label{isot cart biin}
Any isotropic relation $X \overset{f}{\longleftarrow} Y$ can be decomposed as the direct sum of a cartesian  relation $X_C \overset{f_C}{\longleftarrow} Y_C$ and a biinjective relation $X_B \overset{f_B}{\longleftarrow} Y_B$, both of which are isotropic.
\end{prop}

\pf
To begin, we decompose the kernel and indeterminacy of $f$ as $0f = (0f \cap R_Y) \oplus Y_0$ and $ f0 = (f0 \cap R_X) \oplus X_0$,
where $Y_0$ and $X_0$ are arbitrary complements of the respective intersections with the radicals. 

Next, we choose complements $Y_B$ of $0f$ in $(0f)^{\perp}$ and $X_B$ of $f0$ in $(f0)^{\perp}$.  Because $f$ is isotropic, $Xf \subseteq (0f)^{\perp}$ and $fY \subseteq (f0)^{\perp}$, and because $0f \subseteq Xf$ and $f0 \subseteq fY$, we have $Xf = 0f \oplus (Xf\cap Y_B)$ and $fY = f0 \oplus( fY \cap X_B)$ by the modular law.\footnote{The modular law is the fact that, for subspaces $E$, $F$ and $G$ of a vector space $V$, if $E \subseteq G$, then {$G \cap (E + F) = E + (G \cap F)$.} See \cite{ro:advanced}, p. 56., for example.}
Furthermore, 
\begin{equation*}
f0 \times 0f \subseteq f \subseteq Xf \times fY \subseteq [f0 \times 0f] \oplus [X_B \times \ybar_B],
\end{equation*}
so, again by the modular law, $f = [f0 \times 0f] \oplus f \cap [X_B \times \ybar_B]$.
As a complement of $f0 \times 0f$ in $f$, the second summand in this decomposition is a biinjective relation; we denote it by $f_B$.

Finally, we choose a complement of $(0f)^{\perp}$ in $Y$; it is paired non-degenerately with $Y_0$ by the presymplectic structure of $Y$, so we call it $Y_0^*$; similarly, we choose a complement $X_0^*$ of $(f0)^{\perp}$  in $X$. This gives a decomposition
\begin{equation*}
X \times \ybar = [(X_0 \oplus X_0^* \oplus( f0 \cap R_X)) \times \overline{(Y_0 \oplus Y_0^* \oplus (0f \cap R_Y ))}] \oplus [X_B \times \ybar_B].
\end{equation*}
The intersection of $f$ with the first summand in square brackets is simply the cartesian relation $f0 \times 0f$. We denote this part of $f$ by $f_C$ and obtain the desired decomposition $f = f_C \oplus f_B$. 
\qed

\begin{rmk}
\label{cart-decomp}
\emph{
The cartesian isotropic relation $f_C$ in the proof above can itself be written as a direct sum of:
 \bigskip \newline
(1) the following isotropic relations between symplectic spaces
\begin{itemize}
\item $X_0 \oplus X_0^* \leftarrow \zero$ given by  $X_0 \times \zero$, where $X_0$ and $X_0^*$ are lagrangian subspaces of the target space,
\item $\zero \leftarrow Y_0 \oplus Y_0^*$ given by  $\zero \times Y_0$, where $Y_0$ and $Y_0^*$ are lagrangian subspaces of the source space; 
\end{itemize}
(2) the following isotropic relations between spaces with the zero presymplectic structure
\begin{itemize}
\item $f0 \cap R_X \leftarrow \zero$ given by  $(f0 \cap R_X) \times \zero$,
\item $\zero \leftarrow 0f  \cap R_Y$ given by $\zero \times (0f \cap R_Y)$. 
\end{itemize}
}

\emph{A general cartesian relation $X \overset{c}{\longleftarrow} Y$ may not be decomposable in this manner (as $f_C$ is), but any cartesian relation \emph{can} be decomposed as a direct sum of summands of these four types, plus zero relations of the types listed in Proposition \ref{prop-decomp-biin} below. In such a decomposition, one may think of the summands which are zero relations as the `biinjective part' of $c$, \text{i.e.} corresponding to the biinjective summand of the decomposition given in Proposition \ref{isot cart biin}. It straightforward to see that zero relations are the only kinds of relation which are both cartesian and biinjective.
}
\end{rmk}

\section{Decomposition of biinjective relations}

The previous proposition and remark showed that any isotropic relation can be decomposed as the direct sum of four cartesian isotropic relations which have a simple form, plus a fifth isotropic relation which is biinjective. Thus, we now turn to this latter kind of relation. 

\begin{prop}
\label{prop-decomp-biin}
Any biinjective isotropic relation $X \overset{g}{\longleftarrow} Y$ between presymplectic spaces can be decomposed as a direct sum of: 
 \bigskip \newline
(1) the following isotropic relations between symplectic spaces
\begin{itemize}
\item the graph of a symplectomorphism $X_S \overset{g_S}{\from} Y_S,$
\item zero relations $X_S' \from \zero$ and $\zero \from Y_S',$
\item a partially defined biinjective relation $X_L \oplus X_L^* \from Y_L \oplus Y_L^*$ given by the graph of an isomorphism $X_L \overset{g_L}{\from} Y_L$, where $X_L, X_L^*$ and $Y_L, Y_L^*$ are lagrangian subspaces of the target and source spaces, respectively. 

\end{itemize}
(2) the following isotropic relations between zero-presymplectic spaces
\begin{itemize}
\item an isomorphism $X_R \overset{g_R}{\from} Y_R,$
\item zero relations $X_R' \from \zero$ and $\zero \from Y_R'$; and
\end{itemize}
(3) the following isotropic relations where the target space is symplectic and the source space is zero-presymplectic, or vice versa. In each case, the symplectic space has twice the dimension of the corresponding zero-presymplectic space:


\begin{itemize}
\item an injective everywhere-defined relation $gY_I \oplus (gY_I)^* \from  Y_I$ given by the graph of an isomorphism $gY_I \overset{g_{Y_I}}{\longleftarrow} Y_I$, where $gY_I$ and $(gY_I)^*$ are lagrangian subspaces of the target space. 
\item a surjective partially-defined relation $X_I \from X_Ig \oplus (X_Ig)^*$ given by the graph of an isomorphism $X_I  \overset{g_{X_I}}{\longleftarrow} X_Ig$, where $X_Ig$ and $(X_Ig)^*$ are lagrangian subspaces of the source space.  
\end{itemize}
\end{prop}

\pf
In the following, we will often think of the biinjective relation $g$ as an isomorphism to the image $gY$  from the domain $Xg$ of $g$. 

First of all, we decompose the radicals $R_Y$ and $R_X$. Consider the image of the radical $R_Y$ under $g$, \text{i.e.} $g(Xg \cap R_Y)$, and its intersection with the radical in $X$, $X_R := g(Xg \cap R_Y) \cap R_X$. Let $Y_R$ denote the preimage of $X_R$ under $g$. Then $g$ restricts to an isomorphism $g_R$ from $Y_R \subseteq Xg \cap R_Y$ to $X_R \subseteq gY \cap R_X$. We choose subspaces $Y_I$ and $X_I$ such that $ Xg \cap R_Y = Y_R \oplus Y_I$ and $gY \cap R_X = X_R \oplus X_I$,
and denote by $g_{Y_I}$ and $g_{X_I}$ the restrictions of $g$ to $Y_I$ and $X_Ig$, respectively.  Note that $X_Ig \cap R_Y = 0$ and $gY_I \cap R_X = 0$, since $g$ is biinjective. As an isotropic relation, $g$ relates isotropic subspaces to isotropic subspaces; thus $X_Ig$ and $gY_I $ are isotropic because $X_I$ and $Y_I$ are. To complete our decomposition of the radicals, we choose subspaces $Y_R'$ and $X_R'$ such that $R_Y = Y_R \oplus Y_I \oplus Y_R'$ and $R_X = X_R \oplus X_I \oplus X_R'$. 

Next, we construct decompositions of the rest of $Xg$ and $gY$. Any complement to $Y_R  \oplus Y_I \oplus X_Ig$ in $Xg$ is itself a presymplectic space; we choose such a complement and decompose it as a direct sum of its radical $Y_L$ and a complement $Y_S$ to this radical (so $Y_S$ is a symplectic subspace). By biinjectivity, the image of $Y_L \oplus Y_S$ under $g$ is a complement of $X_R \oplus X_I \oplus gY_I$ in $gY$ and has a corresponding induced decomposition $X_L \oplus X_S$. Because $g$ is isotropic, $X_L$  is the radical of $X_L \oplus X_S$ and $X_S$ is symplectic. Note that $gY_I$, $X_L$ and $X_S$ are pairwise orthogonal because their preimages under $g$ are (recall that $g^{-1}(gY_I) \subseteq R_Y$) and the orthogonality of subspaces in the domain and image, respectively, is preserved under isotropic relations. Similarly, $X_Ig$, $Y_L$ and $Y_S$ are pairwise orthogonal. 

We now extend our decompositions to all of $Y$ and $X.$ Let $E_Y$ be a complement of $R_Y$ in $Y$ which contains $X_Ig \oplus Y_L \oplus Y_S$ and let $E_X$ be a complement of $R_X$ in $X$ which contains  $gY_I \oplus X_L \oplus X_S$. Because $E_Y$ is symplectic, there exist independent subspaces $(X_Ig)^*$ and $Y_L^*$ of $E_Y$ which are independent of $X_Ig \oplus Y_L \oplus Y_S$ and which are paired nondegenerately with $X_Ig$ and $Y_L$, respectively, by the presymplectic structure of $Y$. The same holds in $X$ for subspaces $(gY_I)^*$ and $X_L^*$ in $E_X$. The subspace $(X_Ig \oplus (X_Ig)^*) \oplus (Y_L \oplus Y_L^*) \oplus Y_S \subseteq E_Y$ is symplectic; let $Y_S'$ denote its unique orthogonal complement in $E_Y$. This gives a decomposition of $E_Y$ as a direct orthogonal sum of four symplectic subspaces, and in total we obtain
\begin{equation*}
Y = Y_R \oplus Y_R' \oplus Y_I \oplus (X_Ig \oplus (X_Ig)^*) \oplus (Y_L \oplus Y_L^*) \oplus Y_S \oplus Y_S'.
\end{equation*}
The same procedure in $E_X$ gives a corresponding decomposition
$$X = X_R \oplus X_R' \oplus X_I \oplus (gY_I \oplus (gY_I)^*) \oplus (X_L \oplus X_L^*) \oplus X_S \oplus X_S'.$$

From these decompositions of $Y$ and $X$ we can in turn write
\begin{align*}
X \times \ybar  =  &  \ [X_R \times \ybar_R] \oplus [X_R' \times 0 ] \oplus [0 \times \ybar'_R] \\ 
& \oplus [X_I \times \overline{(X_Ig \oplus (X_Ig)^*}) ] \oplus 
[(gY_I \oplus (gY_I)^*) \times \ybar_I ] \\ 
 &  \oplus [(X_L \oplus X_L^*) \times \overline{(Y_L \oplus Y_L^*)}] \oplus [X_S \times \ybar_S] \oplus [X_S' \times 0] \oplus [0 \times \ybar_S']  
\end{align*}
and, by construction, $g$ can be written as a direct sum of its restrictions to each of these summands. These give precisely the nine isotropic relations listed in the statement of the proposition. 
\qed

\section{Indecomposables}

We now prove the main result for isotropic relations between presymplectic spaces and for coisotropic relations between Poisson spaces. We first analyze isotropic relations, which are easier to handle; we then obtain the result for coisotropic relations using duality, thinking of Poisson spaces as the duals of presymplectic spaces.

\begin{thm}
\label{thm-indecomp}
Any indecomposable isotropic relation between presymplectic vector spaces is isomorphic to exactly one of the following thirteen relations.
\bigskip

Relations between symplectic spaces:
\begin{itemize}
\item The identity relation $\mathbb{R}^2 \from \mathbb{R}^2$;
\item The relations $\mathbb{R}^2 \from \zero$ and $\zero \from \mathbb{R}^2$ given by $(\mathbb{R},0)\times \zero$ and $\zero \times (\mathbb{R},0)$ respectively;
\item The zero relations $\mathbb{R}^2 \from \zero$ and $\zero \from \mathbb{R}^2$;
\item The relation $\mathbb{R}^2 \from \mathbb{R}^2$ given in coordinates $(q_1, p_1, q_2, p_2)$ by the equations $q_1 = q_2 = 0$ and $p_1 = p_2$. 
\end{itemize}

Relations between zero-presymplectic spaces: 
\begin{itemize}
\item The identity relation $\mathbb{R} \from \mathbb{R}$; 
\item The relations $\mathbb{R} \from \zero$ and $\zero \from \mathbb{R}$ given by $\mathbb{R} \times \zero$ and $\zero \times \mathbb{R}$ respectively; 
\item The zero relations $\mathbb{R} \from \zero$ and $\zero \from \mathbb{R}$.
\end{itemize}

Relations between a symplectic space and a zero-presymplectic space:
\begin{itemize}
\item The relations $\mathbb{R}^2 \from \mathbb{R}$ and $\mathbb{R} \from \mathbb{R}^2$ given by the natural isomorphisms $(\mathbb{R}, 0) \from \mathbb{R}$ and $\mathbb{R} \from (\mathbb{R},0)$ respectively. 
\end{itemize}
Any indecomposable coisotropic relation between Poisson vector spaces is isomorphic to exactly one of the following thirteen relations.

Relations between non-degenerate Poisson (i.e., symplectic) spaces:
\begin{itemize}
\item The identity relation $\mathbb{R}^2 \from \mathbb{R}^2$;
\item The relations $\zero \from \mathbb{R}^2$ and $\mathbb{R}^2 \from \zero$ given by $\zero \times (\mathbb{R},0)$ and $(\mathbb{R},0)\times \zero$ respectively;
\item The relations $\zero \from \mathbb{R}^2$ and $\mathbb{R}^2 \from \zero$ given by $\zero \times \mathbb{R}^2$ and $\mathbb{R}^2\times \zero$ respectively;
\item The relation $\mathbb{R}^2 \from \mathbb{R}^2$ given in coordinates $(q_1^*, p_1^*, q_2^*, p_2^*)$ by the equation $p_1^* = p_2^*$. 
\end{itemize}

Relations between zero-Poisson spaces: 
\begin{itemize}
\item The identity relation $\mathbb{R} \from \mathbb{R}$; 
\item The zero relations $\zero \from \mathbb{R}$ and $\mathbb{R} \from \zero$;
\item The relations $\zero \from \mathbb{R}$ and $\mathbb{R} \from \zero$ given by $\zero \times \mathbb{R}$ and $\mathbb{R} \times \zero$ respectively.
\end{itemize}

Relations between a non-degenerate Poisson space and a zero-Poisson space:
\begin{itemize}
\item The relations $\mathbb{R} \from \mathbb{R}^2$ and $\mathbb{R}^2 \from \mathbb{R}$  given by the graph of the natural projection onto the first factor $\mathbb{R} \from \mathbb{R}^2$ and its transposed relation. 
\end{itemize}
Any isotropic or coisotropic relation may be decomposed as a direct sum of indecomposable relations. The multiplicities $n_1, ... , n_{13}$ of the indecomposables are invariants of the relation itself. Two (co)isotropic relations are isomorphic if and only if their corresponding multiplicities are equal. 
\end{thm}

\pf
We begin with the case of isotropic relations. First, we verify that each of the listed relations is indeed indecomposable. 

For the first five relations, this follows from the fact that $\reals^2$ is indecomposable as  a symplectic vector space.    Decomposing the sixth relation would require writing it as a direct sum of a
relation $\reals^2 \from \zero$ and a relation $\zero \from \reals^2$, which is impossible, since such a relation would have to be the product of isotropic subspaces in the target and source.

Next, if the identity relation $\mathbb{R} \from \mathbb{R}$ were decomposable, it would have to be the sum of a surjective relation $\mathbb{R} \from \zero$ and an everywhere-defined relation $ \zero \from \mathbb{R}$. The sum of such relations would be $2$-dimensional, but the identity relation is only $1$-dimensional, so a non-trivial decomposition is indeed impossible. The next four relations are clearly also indecomposable for dimension reasons. 

For the last two relations, consider first the relation $\mathbb{R}^2 \from \mathbb{R}$. A non-trivial decomposition would have to consist of a relation $\mathbb{R}^2 \from \zero$ with an isotropic subspace as its image, and an everywhere-defined relation $\zero \from \mathbb{R}$. But then their sum would be a $2$-dimensional relation, whereas the relation in question is $1$-dimensional, a contradiction. The indecomposability of the last relation, $\mathbb{R} \from \mathbb{R}^2$, can be shown analogously.

Next, we prove that any isotropic relation may be decomposed as a direct sum of indecomposable relations of the types listed in the theorem (which shows, in particular, that the list is complete). By Proposition \ref{isot cart biin} and Remark \ref{cart-decomp}, it is enough to show decomposability for the four types of cartesian relations listed in Remark \ref{cart-decomp} and the nine types of biinjective relations listed in Proposition \ref{prop-decomp-biin}. 

The cartesian relations and the zero relations are isomorphic, respectively, to copies of the four cartesian relations and the four zero relations in the statement of the theorem.  The graph of a symplectomorphism can be decomposed into copies of the identity relation on $\mathbb{R}^2$, and  the graph of an isomorphism between lagrangian subspaces of symplectic spaces can be decomposed into copies of the second relation $\mathbb{R}^2 \leftarrow \mathbb{R}^2$ in the statement.  
An isomorphism relation   between zero-presymplectic spaces is isomorphic to a direct sum of copies of the indecomposable identity relation on $\mathbb{R}$.   Finally, the last two types of isotropic relation listed in Proposition \ref{prop-decomp-biin} are isomorphic to multiples of the indecomposable relations $(\mathbb{R},0) \from \zero$ and $\zero \from (\mathbb{R},0)$ respectively.

We now show that, for any decomposition of a given isotropic relation $X \overset{f}{\from} Y$ into indecomposables, the multiplicities of the indecomposable types are the same. 
Following the approach in \cite{lo-we:coisotropic}, we write down a list of invariants associated to $f$: 

\begin{equation*}\label{canonical invariants}
\def\arraystretch{1.1}
\begin{array}{rclclll}
k_1 & = & \tfrac{1}{2}\dim (X/{R_X})  & \quad & k_7&=& \dim (f0 \cap R_X)  \\
k_2 	& = & \tfrac{1}{2}\dim (Y/{R_Y})  & \quad & k_8& = & \dim (0f \cap R_Y)  \\
k_3  	& = & \dim R_X  & \quad & k_9 & = & \dim (fY \cap R_X)  \\
k_4 	& = & \dim R_Y  &\quad & k_{10} & = & \dim (Xf \cap R_Y)  \\
k_5  	& = & \dim f0  &\quad & k_{11} & = & \dim (f \cap (R_X \times R_Y))  \\
k_6  	& = & \dim  0f   & \quad & k_{12} & = & \tfrac{1}{2} \dim  ( Xf/ (Xf \cap (Xf ) ^{\perp})) \\
& & & \quad & k_{13} & = & \dim ( (Xf \cap (Xf)^{\perp})/0f ) 
\end{array}
\end{equation*}

Evaluating these invariants for the indecomposable types, we find that 
the column vector $\bf{k}$ of invariants can be calculated from the column vector $\bf{n}$ of multiplicities by multiplication by the following matrix,

\begin{equation*}
M = 
\left(
\begin{array}{ccccccccccccc}
0& 1& 0& 0& 0& 0& 0& 1& 0& 1& 1& 0& 1\\
1& 0& 0& 0& 0& 0& 0& 0& 1& 1& 1& 1& 0\\
0& 0& 0& 1& 1& 0& 1& 0& 0& 0& 0& 0& 0\\
0& 0& 1& 0& 1& 1& 0& 1& 0& 0& 0& 0& 0\\
0& 1& 0& 1& 0& 0& 0& 0& 0& 0& 0& 0& 0\\
1& 0& 1& 0& 0& 0& 0& 0& 0& 0& 0& 0& 0\\
0& 0& 0& 1& 0& 0& 0& 0& 0& 0& 0& 0& 0\\
0& 0& 1& 0& 0& 0& 0& 0& 0& 0& 0& 0& 0\\
0& 0& 0& 1& 1& 0& 0& 0& 1& 0& 0& 0& 0\\
0& 0& 1& 0& 1& 0& 0& 1& 0& 0& 0& 0& 0\\
0& 0& 1& 1& 1& 0& 0& 0& 0& 0& 0& 0& 0\\
0& 0& 0& 0& 0& 0& 0& 0& 0& 0& 1& 0& 0\\
0& 0& 0& 0& 1& 0& 0& 1& 1& 1& 0& 0& 0
\end{array}
\right)
\
\end{equation*}
which is invertible over $\mathbb{Z}$.
Thus we can recover the multiplicities $\bf{n}$ from the invariants $\bf{k}$ by multiplication by
the inverse matrix $M^{-1}$; i.e., equivalent relations have the same multiplicities. 

This completes the proof for isotropic relations; we turn now to coisotropic relations between Poisson vector spaces. 

When we identify a space $X$ with its double dual $X^{**}$, a presymplectic structure on $X$ becomes equivalent to a (constant) Poisson structure on $X^*$. Isotropic relations $X \overset{f}{\leftarrow} Y$ correspond by duality to coisotropic relations $Y^* \overset{f^*}{\leftarrow} X^*$, where the dual relation $f^*$ is the same as the annihilator $f^{\circ}$ of $f$ with respect to the non-degenerate pairing of $X \times Y$ with  $Y^* \times X^*$ defined by 
$$ \langle (x, y) , (\eta, \xi) \rangle = \xi (x) - \eta (y)$$
(for a more detailed discussion of duality, see \cite{we:(co)isotropic}). 

We have shown above that for any isotropic relation  $X \overset{f}{\leftarrow} Y$, the space $X \times Y$ has a direct sum decomposition into thirteen summands such that each of the intersections of $f$ with these summands is isomorphic to a unique multiple of one of the indecomposable isotropic relations listed in the statement of the theorem. Such a decomposition of $X \times Y$ induces a corresponding decomposition of $Y^* \times X^*$, with respect to which $f^*$ can be decomposed as a direct sum of indecomposables which are dual to the corresponding indecomposable summands of $f$ and which have the same multiplicities as their isotropic counterparts. It is thus sufficient to calculate the duals of the canonical indecomposable isotropic relations, which gives the coisotropic indecomposables listed above. This completes the proof. 
\qed

\begin{rmk}
\emph{We note (as we did similarly in \cite{lo-we:coisotropic}) that the set of possible vectors $\textbf{k}$ is a subset of $\mathbb{Z}_{\geq 0}^{13}$ constrained only by the ten inequalities arising from the condition that $\textbf{n} \in \mathbb{Z}_{\geq 0}^{13}$; the set of possible vectors $\textbf{n}$ of multiplicities is, of course, simply $\mathbb{Z}_{\geq 0}^{13}$}.
\end{rmk}

\end{document}